\input amstex
\documentstyle{amsppt}
\magnification=\magstep1
\hsize=5.2in
\vsize=6.8in
\centerline {\bf ON THE COHOMOLOGY OF ACTIONS OF GROUPS} 
\centerline {\bf BY BERNOULLI SHIFTS}
\vskip .1in
\centerline {\rm by}
\vskip .1in
\centerline {\rm SORIN POPA\footnote"*"{Supported in part by 
NSF-Grant 0100883} and 
ROMAN SASYK}

\address Math Dept  
UCLA, Los Angeles, CA 90095-155505\endaddress
\email popa\@math.ucla.edu, \,
 rsasyk\@math.ucla.edu\endemail

\topmatter

\abstract We prove that if $G$ is a countable, 
discrete group having  
infinite, normal subgroups with the relative 
property (T), then the Bernoulli shift 
action of $G$ on ${\underset 
g \in G \to \Pi} (X_0, \mu_0)_g$, 
for $(X_{0},\mu_{0})$ an arbitrary 
probability space, has first cohomology group isomorphic to the 
character group of $G$.
\endabstract

\endtopmatter

\document
\head 1. Introduction \endhead

Let $G$ be a countable discrete group and $\sigma$ a 
measure-preserving, free, ergodic action of $G$ on a 
probability space $(X, \mu)$.
$\sigma$ induces an action (also denoted by $\sigma$)
of $G$ on the abelian von Neumann algebra 
$A=L^{\infty}(X,\mu)$ 
by $\sigma_g(f):=f\circ \sigma_{g^{-1}}$. 
A common example of such actions are 
the {\it Bernoulli shifts} $\sigma$, 
defined by taking an arbitrary 
probability space $(X_0, \mu_0)$, 
then defining $(X, \mu)={\underset 
g \in G \to \Pi} (X_0, \mu_0)_g$, 
where $(X_0, \mu_0)_g$ are identical copies of 
$(X_0, \mu_0)$, and then letting $\sigma_g$ act on $(X, \mu)$ 
by $\sigma_g((x_h)_h)=(x_{g^{-1}h})_h$. 

A $1$-{\it cocycle} for a 
free, ergodic measure-preserving, 
action $\sigma$ of $G$ on a 
probability space $(X, \mu)$ is a map $w:G\rightarrow \Cal U(A)$ 
satisfying the relations $w_{gh}=w_g \sigma_g(w_h)$, 
$\forall g, h\in G_0$ and $w_e=1$, where 
$\Cal U(A)$ is the group of $\Bbb T$ 
valued functions in $A=L^\infty(X, \mu)$. 
For example, any 
character $\gamma$ of $G$ gives a 1-cocycle for $\sigma$ by $w_g=\gamma(g)1, 
g\in G$. 
A 1-cocycle $w$ is {\it co-boundary} 
if there exists $v\in \Cal U(A)$ such that $w_g = 
v \sigma_g(v^*), \forall g$. 
Denote by $Z^1(\sigma)$ the set of all 1-cocycles 
and by $B^1(\sigma)$ the set of co-boundaries. $Z^1(\sigma)$ is 
clearly a commutative group under multiplication, with $B^1(\sigma)$ 
a subgroup. The corresponding quotient group $H^1(\sigma)= 
Z^1(\sigma)/B^1(\sigma)$ 
is called the {\it first cohomology group} of $\sigma$, and is clearly 
a conjugacy invariant for $\sigma$. 

In the early 80's Klaus Schmidt proved that the group $G$ has the 
property (T) of Kazhdan ([K]) if and only if 
$H^1(\sigma)$ is countable for any free, ergodic, measure 
preserving action $\sigma$ of $G$ ([S2]). He also showed that 
$G$ is amenable iff the Bernoulli shift 
actions of $G$  are non-strongly ergodic, and iff all measure-preserving 
actions of the group $G$ are non-strongly ergodic. 
Related to these results, Connes and Weiss 
proved that $G$ has the property (T) iff all its ergodic, 
free measure-preserving actions are strongly ergodic.

In this paper we obtain the first actual computations 
of cohomology groups $H^1(\sigma)$, in the case 
$\sigma$ is a Bernoulli shift action and the group $G$ 
is {\it weakly rigid} in the following sense:  
$G$ contains  
infinite, normal subgroups $H \subset G$ such that 
$(G, H)$ has the relative property (T) of 
Kazhdan-Margulis([M], [dHV]), 
i.e., any representation of $G$ that weakly contains 
the trivial representation of $G$  must contain 
the trivial representation of 
$H$. Note that any group $G$ of the form $G=H \times \Gamma$ 
with $H$ an infinite group with the property (T) of Kazhdan 
is weakly rigid.

\proclaim{Theorem} If $G$ is a countable, weakly rigid
discrete group 
and   
$\sigma$ is a Bernoulli shift action of $G$ then $H^1(\sigma)$ 
is equal to the character group of $G$. 
\endproclaim

\proclaim{Corollary} If $\Gamma$ is an arbitrary countable 
discrete abelian group, $G=SL(n, \Bbb Z) \times \Gamma$, 
for some $n \geq 3$, and 
$\sigma$ is a Bernoulli shift action of $G$ then 
$H^1(\sigma)=\hat{\Gamma}$. 
\endproclaim

We mention that the similar result for (purely) non-commutative 
Bernoulli shifts was obtained in ([Po]). 
In fact, to prove the above Theorem we  
will follow the line of arguments in ([Po]), with the   
commutativity  allowing 
many simplifications.

\vskip .3in

\head 2. Preliminaries \endhead
\vskip .1in

Let $G$ be a discrete group and $\sigma$ a measure preserving 
action of $G$ on a standard probability measure space 
$(X,\mu)$. The action it implements on the 
abelian von Neumann algebra $A=L^\infty(X,\mu)$, 
still denoted by $\sigma$, preserves the integral and thus 
extends to an action (or unitary representation) 
$\sigma$ of $G$ on the Hilbert space $L^2(X, \mu)$. 
We denote by $\Cal U(A)$ the group of unitary elements 
in $A$. Besides the notion of
1-cocycles for $\sigma$  defined in the introduction
we need the following:
\vskip .1in
\noindent
{\bf 2.1. Definition.} A {\it weak 1-cocycle} for the 
action $\sigma$ is a function 
$w :G\rightarrow \Cal U(A)$ satisfying  	
$w_{gh}=w_g\sigma_g(w_h)\mod \Bbb C, \forall g, h\in G,$ and $w_e=1$.
A weak cocycle $w$ is a {\it weak coboundary} 
if there exist a unitary $u$ in $A$ 
such that $w_g=u\sigma_g(u^*)	\mod \Bbb C, \forall g\in G$.
Note that if $w$ is a weak 1-cocycle for $\sigma$ and 
$v\in \Cal U(A)$ then 
$w'_g=vw_g\sigma_g(v^*), g\in G$ is also a weak 
1-cocycle for $\sigma$. Two weak 1-cocycles $w, w'$ for which  
there exists $v$ as above are called {\it equivalent}. 

\vskip .1in
\noindent
{\it 2.2. Remarks}. 1$^\circ$. Let $w$ be a weak 1-cocycle for $\sigma$ 
and denote by $\gamma(g,h)\in \Bbb T$ the scalar satisfying  
$w_{gh}=\gamma(g,h)w_g\sigma_g(w_h), \forall g,h$. Condition  
$w_e=1$ then implies 
$\gamma(e,g)=\gamma(g,e)=1,  \forall g\in G$.
Also, the associativity relation 
$w_g(w_hw_k)=(w_gw_h)w_k$ entails   
$$
\gamma(g,h)\gamma(gh,k)=\gamma(g,hk)\gamma(h,k), \forall g,h,k. 
$$
A function $\gamma:G\times G\rightarrow\Bbb T$ 
that verifies the previous conditions 
is called a {\it scalar valued} (or $\Bbb T$-valued) 
{\it 2-cocycle} for the group $G$. Thus  any 
weak 1-cocycle $w$ for the action $\sigma$ 
has associated a scalar 2-cocycle $\gamma=\gamma_w$.

2$^\circ$. If the weak 1-cocycle $w$ is a weak coboundary and  
$w_g=\lambda_g v\sigma_g (v^*)$ 
with $\lambda_g \in \Bbb C$ then 
$\lambda_{gh}=\gamma(g,h)\lambda_g\lambda_h$. 
In particular if $w$ is a genuine 1-cocycle then
$\lambda$ follows a character of $G$.

\proclaim{2.3. Lemma}
Let $w$  be a weak 1-cocycle for the action $\sigma$,   
with scalar 2-cocycle $\gamma$.

$1^\circ$. For $g \in G$ and $\xi \in L^2(X,\mu )$ denote 
$\sigma_g^w(\xi):=w^*_g \sigma_g(\xi).$
Then $\sigma^w$ is a projective representation 
of $G$ on $L^2(X,\mu)$ with scalar 
$2$-cocycle $\overline{\gamma}$.

$2^\circ.$ Let $\Cal{HS}$ the space of Hilbert-Schmidt 
operators on $L^2(X,\mu)$ and 
for each $T$ in $\Cal{HS}$ denote  
$\widetilde{\sigma}_g^w(T):=\sigma_g^wT\sigma_g^{w *}.$
Then $\widetilde{\sigma}^w$ is a unitary representation of 
$G$ on $\Cal{HS}.$

$3^\circ$. If we identify an element $T\in \Cal{HS}$ 
with an element $a$ of $L^2(X, \mu)
\overline{\otimes} L^2(X, \mu) 
\simeq L^2(X \times X, \mu \times \mu)$ 
in the usual way, 
then $\widetilde{\sigma}^w_g(T)=1\otimes w_g
\cdot (\sigma_g \otimes\sigma_g)(a)
\cdot w_g^*\otimes 1$.

\endproclaim

\noindent
{\it Proof}. $1^\circ$. We have:

$$
\sigma_g^w\sigma_h^w(\xi)
=\sigma_g^w\sigma_h(\xi)w_h^*
=\sigma_g(\sigma_h(\xi))\sigma_g(w_h^*)w_g^*
$$
$$
=\sigma_{gh}(\xi)\overline{\gamma(g,h)}w_{gh}^*
=\overline{\gamma(g,h)}\sigma_{gh}^w(\xi), 
$$ 
showing that $\sigma^w$
is a projective representation.

$2^\circ$. If $T\in \Cal{HS}$ then  
$$
\widetilde{\sigma}_g^w\widetilde{\sigma}_h^w(T)
=\sigma_g^w\sigma_h^wT\sigma_h^{w *}\sigma_g^{w *}
=\overline{\gamma(g,h)}\sigma_{gh}^w T\gamma(g,h)\sigma_{gh}^{w *}
=\widetilde{\sigma}_{gh}^w(T).
$$

$3^\circ$. Since the space $\Cal{HS}$
of Hilbert-Schmidt operators on $L^2(X, \mu)$ 
is isomorphic to $L^2(X, \mu)\overline{\otimes} 
L^2(X, \mu)$ 
via the identification  
$x\otimes y^{*}(\xi)= (\xi|x)\,y$, we have: 

$$
\widetilde{\sigma}^w_g(x\otimes y^*)(\xi)=
(\sigma_g^{w^*}(\xi)|x)\sigma_g^w(y)=
(\xi|{\sigma}_g^w(x))\sigma_g(y)w_g^*
$$
$$
=(\sigma_g(x)w_g^*)\otimes(w_g\sigma_g(y)^*)(\xi)
=(1\otimes w_g)(\sigma_g(x)\otimes
\sigma_g(y^*))(w_g^*\otimes 1)(\xi).
$$
Then extend by linearity.
\hfill $\square$

\proclaim{2.4. Lemma} With the notations 
of ${\text{\rm Lemma 2.3}}$, the following are equivalent:

$(1).$ $\widetilde{\sigma}^w$ contains 
a copy of the trivial representation.

$(2).$ $\sigma^w$ has a non trivial, 
invariant finite dimensional subspace 
$\Cal H_0 \subset L^2(X,\mu )$. 

$(3).$ There exist  $a\ne 0$ in $L^2(X)
\overline{\otimes} L^2(X)$
such that $(1\otimes w_g)\cdot (\sigma_g\otimes
\sigma_g)(a)\cdot (w_g^*\otimes 1)=a$.

\endproclaim

\noindent
{\it Proof}. (1)$\Rightarrow$ (2). First note that the action 
$\widetilde{\sigma}^w$ can be extended to all
$\Cal B (L^2(X,\mu))$.
Also note that if $T\in \Cal {HS}$ is 
fixed by $\widetilde{\sigma}^w$, so is $T^*$.
Thus the trace class operator
$TT^*\in \Cal B (L^2 (X,\mu))$ 
is fixed by $\widetilde{\sigma}^w$.
By the Borel functional calculus, all the spectral
projections of $TT^*$ are fixed by 
$\widetilde{\sigma}^w$. As they have finite trace,
it follows that they are projections on finite 
dimensional subspaces.
Choose $\Cal H_0$ a non-trivial finite
dimensional subspace of $L^2(X,\mu)$ 
corresponding to some spectral projection 
$P$ of $TT^*$, and note that 
$\widetilde\sigma_g ^w(P)$ is the projection of
$L^2(X,\mu)$ onto $\sigma_g ^w(\Cal H_0)$. Then as
$\widetilde{\sigma}_g^w(P)=P$, it follows that
${\sigma}_g^w(\Cal H_0)=\Cal H_0.$ Thus
$\Cal H_0$ is an invariant subspace for $\sigma^w$.

(2)$\Rightarrow$ (1). If $\Cal H_0$
is an invariant subspace for $\sigma^w$, take $P$ the 
finite rank projection of $L^2(X,\mu)$ onto $\Cal H_0$.
Then $P$ is invariant for the action $\widetilde\sigma^w$.

(2)$\Leftrightarrow$ (3). Trivial by previous lemma.
\hfill{\hfill $\square$}

Recall that the measure preserving action $\sigma$ 
of $G$ on the probability space $(X, \mu)$ 
is weakly mixing if and only if 
the only finite dimensional subspace 
of $L^2(X, \mu)$ invariant to $\sigma$ is $\Bbb C1$ 
(see e.g. [BMe]).

\proclaim{2.5. Lemma} Assume $\sigma$ 
is weakly mixing and let $w$ be a weak $1$-cocycle for $\sigma$. 
Then $\tilde{\sigma}^w$ contains a copy of the 
trivial representation if and only if $w$ is a weak coboundary. 
Moreover, if this is the case, then the unitary 
element $u\in L^{\infty}(X, \mu)$ with 
$w_g = u\sigma_g(u^*)$ $\mod \Bbb C$, 
$\forall g\in G$, 
is unique up to a scalar multiple.
\endproclaim
\noindent
{\it Proof}. If $w_g = \lambda_g u^*\sigma_g(u), \forall 
g\in G,$ then $\sigma^w_g(\Bbb Cu)= \Bbb C \sigma_g(u)w_g^* 
= \Bbb Cu$, thus $\Cal H_0 = \Bbb Cu$ is $\sigma^w$-invariant.  
Thus, by Lemma 2.4, the orthogonal 
projection onto $\Cal H_0$ is a 
fixed point for $\tilde{\sigma}^w$. 

Conversely, let $\Cal H_0\subset L^2(X, \mu)$ 
be a $\sigma^w$-invariant finite dimensional subspace.  
Choose an orthonormal basis 
$\{ \xi_1,\ldots ,\xi_n\}$ of $\Cal H_0$ and note that 
$\{\eta_i\}_i = \{\sigma_g(\xi_i)w_g^*\}_i$ is also an 
orthonormal basis of $\Cal H_0$. But an easy 
computation shows that  
$\Sigma_i \xi_i \xi_i^* = \Sigma_i \eta_i\eta_i^* \in L^1(X, \mu)$ 
for any two orthonormal basis of $\Cal H_0$. 
Thus, since 
$$
\sigma_g(\Sigma_i \xi_i\xi_i) = 
\Sigma_i 
(\sigma_g(\xi_i)w_g^*)(w_g \sigma_g(\xi^*))
$$
$$
= \Sigma_i \eta_i \eta_i^* =\Sigma_i \xi_i\xi_i^*
$$ 
and since $\sigma$ is ergodic on $(X, \mu)$, it follows that 
$\Sigma_i \xi_i \xi_i^* \in \Bbb C1$. In particular, 
all $\xi_i$ are bounded elements, $\xi_i \in L^\infty(X, \mu)$. 

But since 
$$ 
\sigma_g(\xi_i\xi_j ^* )=
\sigma_g(\xi_i)w_g^* w_g\sigma_g(\xi_j ^*)
=\sigma_g^w(\xi_i)(\sigma_g^w(\xi_j))^*, 
$$
the finite dimensional subspace $\Cal H_0 \cdot \Cal H_0^*$
of $L^{\infty} (X,\mu)$ spanned by
$\{ \xi_i \xi_j^*\}_{i,j=1}^n$ 
is $\sigma$-invariant. Since $\sigma$ 
is weakly mixing, it follows that $\Cal H_0\Cal H_0^*=\Bbb C1$. 
Thus $\xi_i\xi_j^*\in \Bbb C1, \forall i,j,$ 
which implies that $n=1$ and $\xi_1$
is a scalar multiple of a unitary element.

Finally, if $w_g=u^*\sigma_g(u) \mod\Bbb C$ and 
$w_g={u'}^*\sigma_g(u^{'})\mod\Bbb C$,
then $u^{'}u^*= \sigma_g(u^{'}u^*)\mod \Bbb C$, 
i.e. the subspace $\Bbb Cu^{'}u^*$ 
is invariant to $\sigma$,
implying that $u'=u \mod{\Bbb T}$.
\hfill $\square$

\vskip .3in

\head 3. The main result \endhead

Let $(X_0,\mu_0)$ be a nontrivial probability space and $G$
an infinite discrete group. Denote 
$(X,\mu):=\prod_{g \in G}(X_0,\mu_0)_g$ and let $\sigma$ 
be the action of $G$ on 
$(X,\mu)$ by $G$-{\it Bernoulli shifts}, i.e., 
$\sigma_g((x_h)_h)= (x_{g^{-1}h})_h$. 
This action is well known to be mixing. Note that if $H \subset G$ 
is a subgroup of $G$ then $\sigma_{|H}$ is a $H$-Bernoulli 
shift. Also, note that the (diagonal) product of 
two $G$-Bernoulli shifts is a $G$-Bernoulli shift. 

Recall from ([dHV]) that 
an inclusion of discrete groups $H \subset G$ has the 
{\it relative property} $(T)$ if the following condition 
holds true: 

\vskip .1in
\noindent
{\it 3.0.} There exist a finite set of elements
$g_1,g_2,\dots ,g_n$ in $G$ and $\epsilon > 0$ such that
for any unitary representation $\pi$ of $G$
on a Hilbert space $\Cal H$ 
which has a unit vector $\xi$ with  
$\| \pi(g_i)\xi-\xi\|_{\Cal H}
<\epsilon$ for all $1\le i\le n$, 
there exists a unit vector fixed by $\pi_{| H}$.

\vskip .1in

By a result of Jolissaint ([Jo]), the above condition  
is equivalent to the following: 
\vskip .1in
\noindent
{\it 3.0'}. Given any $\epsilon>0$ there exist a finite 
set of elements $g_1,g_2,\dots ,g_n$ in $G$ and  
$\delta > 0$ such that 
for any unitary representation $\pi$
of $G$ on a Hilbert space $\Cal H$ 
that has a unit vector $\xi$ such that
$$
\| \pi(g_i)\xi-\xi\|_{\Cal H}
<\delta \text{ for all }  1\leq i\leq n
$$ 
then
$$ 
\| \pi(g)\xi-\xi\|_{\Cal H}
<\epsilon \text{ for all }  g\in H.
$$

\vskip .1in
\proclaim{3.1. Theorem} Let $G$ be a countable discrete 
group and $H \subset G$ a subgroup 
with the relative property $(T)$. Given any
weak 1-cocycle $w$ for a $G$-Bernoulli
shift $\sigma$, $w|_H$ is a weak coboundary.
\endproclaim

\noindent
{\it Proof}. We first prove the case when 
$(X_0,\mu_0)$ is non atomic, thus isomorphic to 
$(\Bbb T,\lambda)$, the torus with its Haar measure. 

Denote by $A$ the abelian von Neumann Algebra 
$ L^\infty(X, \mu).$ By Lemma 2.5, 
it is sufficient to 
prove that there exists 
$u \in \Cal U(A\overline{\otimes} A)$ 
such that $\tilde{\sigma}^w_h(u)=u, \forall h\in H$. We'll 
prove this in the Lemmas 3.2-3.5 below. 

\proclaim{3.2. Lemma}. There exists a continuous action 
$\alpha$ of $\Bbb R$ on $A\overline{\otimes} A 
\simeq L^\infty(X\times X, \mu \times \mu)$, by automorphisms 
preserving the integral over $\nu \times \nu$, such that:
\vskip .05in
\noindent
$(3.2.1)$. $\alpha$ commutes with the Bernoulli shift 
$\tilde{\sigma}=\sigma\otimes\sigma.$
\vskip .05in
\noindent
$(3.2.1)$. $\alpha_1(A\otimes\Bbb C)=\Bbb C\otimes A.$
\endproclaim
\noindent
{\it Proof}. Denote $A_0 = L^\infty(\Bbb T, \lambda)$, 
$\tilde{A}_0=A_0  
\overline{\otimes} A_0$ and $\tau_0$ 
the functional on $\tilde{A}_0$ given by the integral 
over $\lambda \times \lambda$. 
We first construct a continuous action 
$\beta : \Bbb R \rightarrow {\text{\rm Aut}} (\tilde{A}_0, \tau_0)$ 
such that $\beta_1(A_0 \otimes \Bbb C)=\Bbb C \otimes A_0$.  

Let $u$ (resp. $v$) be a Haar unitary generating   
$A_0 \otimes \Bbb C \simeq L^\infty(\Bbb T,\lambda)$ (resp. 
$\Bbb C \otimes A_0$). 
Thus, $u, v$ is a pair of generating Haar unitaries 
for $\tilde{A}_0$, i.e., 
$\{u^nv^m\}_{n, m \in \Bbb Z}$ is an orthonormal basis for $L^2(\tilde{A}_0, \tau_0) \simeq L^2(\Bbb T, \lambda) \overline{\otimes}  
L^2(\Bbb T, \lambda)$. We need to construct 
the action $\beta$ so that $\beta_1(u)=v$. 

Note that given any other 
pair of generating Haar unitaries $u', v'$ for $\tilde{A}_0$, 
the map $u\mapsto u, v\mapsto v'$ extends to a 
$\tau_0$-preserving automorphism of $\tilde{A}_0$. Also, note that     
$v, uv$ is a pair 
of generating Haar unitaries for $\tilde{A}_0$. 
Thus, in order to get $\beta$, it is sufficient to 
find a  continuous action 
$\beta' : \Bbb R \rightarrow {\text{\rm Aut}} (\tilde{A}_0, \tau_0)$ 
such that $\beta'_1(v)=uv$.  
 
Let $h \in \tilde{A}_0$ be a self-adjoint element 
such that $exp(2\pi i h)= u$. It is easy to see that  
for each $t$, $u$ and  
$exp(2\pi i th)v$ is a pair of Haar unitaries. 
Denote by $\beta'_t$ the automorphism 
$u\mapsto u, 
v\mapsto exp(2\pi i th)v$. We then clearly have 
$\beta'_t \beta'_s = \beta'_{t+s}$, $\forall t, s \in \Bbb R$  
and $\beta'_1(v)=uv$.  

Finally, we take $\alpha$ to be the product 
action $\alpha_t=\bigotimes_{g\in G} (\beta_t)_g, t\in \Bbb R$. 
Since $\alpha$ acts identically on the components of the 
product of the $G$-shifts, it commutes with $\sigma$. 
Also, $\alpha_1$ flips $A\otimes \Bbb C$ onto $\Bbb C \otimes A$ 
because each $(\beta_1)_g$  takes $(A_0)_g \otimes \Bbb C$ 
onto $\Bbb C \otimes (A_0)_g$. 
\hfill $\square$

For the next lemma, note that if $K$ is a convex subset of the 
von Neumann algebra $A\overline{\otimes} A = 
L^\infty(X\times X, \mu \times \mu)$ which is bounded 
in the norm $\|\quad \|=\|\quad \|_\infty$, 
then its closure $\Cal K$ in the $w$-operator 
topology on $A\overline{\otimes} A$ coincides with its closure 
in the norm $\|\quad \|_2$ on $L^2(X\times X, \mu \times \mu)$ 
(with $A\overline{\otimes} A\supset K$ 
regarded as a subset of this Hilbert space).  

\proclaim{3.3. Lemma}
For each $t\in \Bbb R$ let $x_t$ 
be the (unique) element of minimal norm-2 in 
$\Cal K_t:=\overline{co}^{ \|\, \|_2}\{(w_h\otimes 1) 
\alpha_t(w_{h}^*\otimes 1)\}_{h\in H}.$
Then $x_t\in A\overline{\otimes} A$ 
and it satisfies the following conditions:

$1^\circ$. $(w_h\otimes 1)\tilde{\sigma}_h(x_t) = 
x_t \alpha_t(w_h\otimes 1), \forall h\in H$.

$2^\circ$. $x_tx_{t}^*\in \Bbb C\otimes\Bbb C.$
\endproclaim

\noindent
{\it Proof}. $1^\circ$. Since $w_h \sigma_h(w_k) = w_{hk}$, 
mod $\Bbb C$, and the actions $\tilde{\sigma}$, $\alpha$ commute, 
it follows that for all $h, k \in G$ we have 
$$
(w_k\otimes 1) \tilde{\sigma}_k ((w_h\otimes 1) \alpha_t(w_h^*\otimes 1)) 
\alpha_t(w_k^*\otimes 1) = (w_{kh}\otimes 1) \alpha_t(w_{kh}^*\otimes 1)
$$
showing that for each fixed $k \in H$ 
the unitary operator 
on $L^2(X \times X, \mu \times \mu)=L^2(A\overline{\otimes} A)$ 
given by $x \mapsto 
(w_k \otimes 1) \tilde{\sigma}_k(x) \alpha_t(w_k^*\otimes 1)$ 
takes $\Cal K_t$ into itself. Thus, by the uniqueness 
of the element of minimal norm $\|\quad \|_2$ in $\Cal K_t$, 
it follows that $x_t= 
(w_k \otimes 1) \tilde{\sigma}_k(x_t) \alpha_t(w_k^*\otimes 1)$, 
$\forall k \in H$. 

$2^\circ$. From the proof of $1^\circ$ and the commutativity 
of $A\overline{\otimes} A$ it follows  that for $k \in H$ we have 
$$
\tilde{\sigma}_k(x_tx_t^*) 
= (w_k \otimes 1) \tilde{\sigma}_k(x_tx_t^*) (w_k^*\otimes 1) 
= x_tx_t^*. 
$$
But since $\sigma |_H$ is weakly mixing, $\tilde{\sigma}_{|H}$ 
is ergodic and thus $x_tx_t^*$ follows a scalar.
\hfill $\square$

\proclaim{3.4. Lemma}
Assume $(G,H)$ has the
relative property $(T)$. If $x_t$ are defined as in ${\text{\rm Lemma 
3.3}}$, then   
there exists $t_0 >0$ such that 
$x_t\ne 0$ and  
$u_t=x_t/\|x_t\|$ is a unitary element in $A\overline{\otimes} A$ 
for all $t \in [0,t_0]$.
\endproclaim

\noindent
{\it Proof}. Let $\epsilon>0$. Let 
$g_1,\ldots ,g_n\in G$ and $\delta>0 $ be given by condition 
$(3.0')$. By the continuity of the action $\alpha_t$, 
there exists $t_0>0$ such that if $0< t \leq t_0$ 
then  
$$
\|(w_{g_i}\otimes 1)
\alpha_t(w_{g_i}^*\otimes 1)-1\|_2<\delta, \forall i .  
$$

Fix $t \in (0, t_0]$. Since the action $\tilde{\sigma}$ 
commutes with the automorphism $\alpha_t$, 
it follows that $\tilde{\sigma}_g \times \alpha_t^n$ 
implements an action of $G \times \Bbb Z$ on $A \overline{\otimes} A$ 
which preserves the functional $\tau$ given by the 
integral over $\mu \times \mu$. 

Let  
$N=(A\overline{\otimes} A)\rtimes (G \times \Bbb Z)$ 
be the corresponding group measure space von Neumann 
algebra ([MvN]) with its canonical trace still denoted by $\tau$. 
Let $U_g \in N$ be the canonical unitaries implementing 
$\tilde{\sigma}_g, g\in G,$ and $U_t$ the unitary implementing 
$\alpha_t$. Denote $U'_g = (w_g \otimes 1) U_g, g\in G$ and 
note that $w$ weak cocycle implies 
$U'_g U'_h = U'_{gh}$ mod $\Bbb C$, $\forall g, h \in G$. 

Let $L^2(N, \tau)$ be the Hilbert space obtained 
by completing  $N$ in the Hilbert norm 
$\|x\|_2 = \tau(x^*x)^{1/2}, x\in N,$ 
and consider the representation $\pi$ of $G$ on $L^2(N, \tau)$ 
given by  
$\pi (g)\xi=U_g'\xi {U'_g}^*$. We have 
$$
\|\pi(g)U_t-U_t\|_2=\|U_g'U_tU_g^{'*}U_t^*-1\|_2
$$
$$
=\|(w_g\otimes 1)U_gU_tU_g^*(w_g^*\otimes 1)U_t-1\|_2=
\|(w_g\otimes 1)\alpha_t(w_g^*\otimes 1)-1\|_2.
$$

Taking $g = g_i, i=1,2,..., n$, 
condition $(3.0')$ on $(G,H)$ implies 
$$
\|(w_h\otimes 1)\alpha_t(w_h^*\otimes 1)-1\|_2
=\|\pi(g_i)U_t-U_t\|_2<\epsilon, \forall h\in H 
$$ 
which in turn implies $\|x-1\|_2 \leq \epsilon, \forall x\in \Cal K_t$. 
In particular, $x_t$ satisfies     
$\|x_t-1\|_2<\epsilon$. Thus  $x_t\ne 0$ and by $3.3.2^\circ$ 
the last part follows. 
\hfill $\square$

\proclaim{3.5. Lemma}
There exists $u\in \Cal U(A\overline{\otimes} A)$ such that 
$$
\tilde{\sigma}^w_h(u)=u, \forall h\in H. 
$$
\endproclaim
\noindent
{\it Proof}.
Choose $n\in \Bbb N$ such that 
$1/n<t_0$, where $t_0$ is by 3.4. With   
$u_t$ defined as in Lemma 3.4,  we let  
$u=u_{1/n} \alpha_{1/n}(u_{1/n})\cdots 
\alpha_{1/n}^{n-1}(u_{1/n})$. By 3.3.1$^\circ$ we have 
$(w_h\otimes 1)\tilde{\sigma}_h(u_{1/n}) = 
u_{1/n} \alpha_{1/n}(w_h\otimes 1)$, which by applying 
on both sides $(\alpha_{1/n})^k = \alpha_{k/n}$, $k = 1, 2, 
..., n-1$, gives 
$$
\alpha_{k/n}(w_h\otimes 1)\tilde{\sigma}_h
(\alpha_{k/n}(u_{1/n})) = 
\alpha_{k/n}(u_{1/n}) \alpha_{(k+1)/n}(w_h\otimes 1)
$$
By applying this repeatedly to $u$, we  get  
$$
(w_h \otimes 1) \tilde{\sigma}_h(u) 
= u \alpha_1(w_h \otimes 1) = u (1 \otimes w_h), \forall h\in H, 
$$
or equivalently $\tilde{\sigma}_h^w (u)=u, \forall h\in H$. 
\hfill $\square$

\vskip .1in

This ends the proof of the nonatomic case.
For the atomic case we need the following:

\proclaim{Lemma 3.6}
Suppose $(X_0,\mu_0)$ is an atomic probability space.  
There exists an embedding 
of $L^\infty(X_0, \mu_0)$ into $L^\infty(\Bbb T, \lambda)$ 
with a sequence of diffuse von Neumann  
subalgebras 
$(B_n)_{n\in\Bbb N}$ of $L^\infty(\Bbb T,\lambda)$ 
such that $B_{n+1}\subset B_n$ 
and $L^\infty(X_0,\mu_0)=\bigcap_{n\in \Bbb N}B_n$.
\endproclaim

\noindent
{\it Proof}. Identify $L^\infty(\Bbb T, \lambda)$ 
with $\overline{\bigotimes}_{n \geq 0} L^\infty(X_n, \mu_n)$, 
where $(X_n, \mu_n)=(X_0, \mu_0), \forall n \geq 0$. Also, 
identify the initial algebra $L^\infty(X_0, \mu_0)$ 
with $L^\infty(X_0, \mu_0) \bigotimes_1^\infty  \Bbb C1 
\subset L^\infty(\Bbb T, \lambda)$ and put 

$$
B_n = L^\infty(X_0, \mu_0) (\Bbb C \bigotimes_1^n 1) 
\bigotimes_{j=n+1}^\infty L^\infty (X_j, \mu_j).
$$

Then $B_n$ are clearly diffuse and 
$\cap_n B_n = L^\infty(X_0, \mu_0)$. 
\hfill $\square$

\vskip .1in

With $L^\infty(X_0, \mu_0) \subset B_n \subset L^\infty(\Bbb T, \lambda)$ 
as in Lemma 3.6, 
denote $A=\overline{\otimes}_{g\in G} L^\infty(\Bbb T, \lambda)_g$  
with its 
subalgebras 
$A_0 = \overline{\otimes}_{g\in G} L^\infty(X_0, \mu_0)_g$ 
and $A_n = \overline{\otimes}_{g\in G} (B_n)_g$, $n \geq 1$.  

The $G$-Bernoulli shift $\sigma$ 
on $A_0$ extends  to $G$-Bernoulli shifts 
on $A$ and $A_n$, $n \geq 1$, still denoted $\sigma$. 
If $w: G \rightarrow \Cal U(A_0)$ is a weak 1-cocycle for $\sigma$ as a  
$G$-Bernoulli shift action on $A_0$, then  
$w$ can also be regarded as a 
weak 1-cocycle for the $G$-Bernoulli 
shift action on $A_n, n\geq 1$. 
The non atomic case of Theorem 3.1 implies 
that $w|_H$ is a weak coboundary for $\sigma_{|H}$ 
as an action on  $A_n$. Thus, for each $n \geq 1$  
there exists a unitary element 
$u_n\in A_n$ such that $w_h=
u_n\sigma_h(u_n^*) \mod{\Bbb C}$. 
By Lemma 2.5, $u_n$ is 
unique up to a scalar multiple. 
Since $A_{n+1}\subset A_n$ and
 $\bigcap_{n\in\Bbb N} A_n=A_0$, it follows  
that $\Bbb Cu_n=\Bbb Cu_{n+1}$ 
and finally $\Bbb Cu_n\in A_0$ for all $n\geq 1$. 
Thus $w|_H$ is a weak 1-cocycle for the 
action $\sigma_{|H}$ on $A_0$.
\hfill \hfill $\square$

\head 4. Applications\endhead

As in  the introduction, 
a group $G$ 
is called {\it weakly rigid} if
it contains  
infinite, normal subgroups $H \subset G$ such that 
the pair 
$(G, H)$ has the relative property $(T)$.

\proclaim{4.1. Theorem}
If $G$ is a weakly rigid group then any
weak 1-cocycle for a $G$-Bernoulli
shift is a weak coboundary.
\endproclaim
\noindent
{\it Proof.} By hypothesis, there exists an infinite  
normal subgroup  $H\subset G$ such that 
$(G,H)$ has the relative property $(T)$.
If $w$ is a weak 1-cocycle for the $G$-Bernoulli shift
$\sigma$, then by Theorem 3.1 
there exists $v\in \Cal U(A)$ such that 
$w_h=v\sigma_h(v^*)$, $\mod \Bbb C$, $\forall h \in H$.

Let $w_g'=v^*w_g\sigma_g(v)$.
Then $w '$ is a 
weak 1-cocycle for $\sigma$ and it satisfies $w_h \in \Bbb T1, 
\forall h \in H$. 

For $a\in A$, denote by 
$L_a \in \Cal B(L^2(X,\mu))$ the (left) multiplication 
operator given by 
$L_a(\xi)=a\xi,  \forall \xi\in L^2(X,\mu)$. Then we have 
$$
L_{w'_g}\sigma_g L_{w'_h}\sigma_h(\xi)=
w'_g\sigma_g(w'_h)\sigma_{gh}(\xi)
=w'_{gh}\sigma_{gh}(\xi) \mod \Bbb T. 
$$
Thus  
$$
(L_{w'_g}\sigma_g) (L_{w'_h}\sigma_h)
=L_{w'_{gh}}\sigma_{gh} \mod \Bbb T.
$$
Similarly   
$$
(L_{w'_g}\sigma_g)^* 
=L_{w'_{g^{-1}}} 
\sigma_{g^{-1}} \mod \Bbb T.
$$
This implies  
$$
(L_{w_g'}\sigma_g) (L_{w_h'}\sigma_h) (L_{w_g'}\sigma_g)^*
=w_{ghg^{-1}}'\sigma_{ghg^{-1}}\mod \Bbb T, 
$$ 
for all $g, h \in G$. Since $w'_h$ are scalars for $h \in H$  
and $ghg^{-1} \in H, \forall g$,  
this further implies   
$$
L_{w'_g} \sigma_{ghg^{-1}} L_{{w'_g}^*} = (L_{w_g'}\sigma_g) \sigma_h  (L_{w_g'}\sigma_g)^*
= \sigma_{ghg^{-1}} \quad \mod \Bbb T
$$
Substituting $h$ for $ghg^{-1}$ and applying the first and last 
term of these equalities to the element $\xi=w'_g 
\in L^2(X, \mu)$, it follows that 
$\sigma_h(w'_g) \in \Bbb Cw'_g, \forall h\in H, g\in G$. 
Since the action $\sigma_{|H}$ is weakly mixing, 
it follows that $\Bbb Cw_g'=\Bbb C1$ for all $g\in G$. 
Thus $w_g = v\sigma_g(v^*), \mod \Bbb T$, $\forall g\in G$, 
i.e., $w$  is a weak coboundary.
\hfill $\square$

\proclaim{4.2. Corollary}
Under the same assumptions as in Theorem 4.1, if 
$w$ is a genuine 1-cocycle 
then $w$ is equivalent to a character of $G$ and 
different characters give non equivalent 1-cocycles. 
In other words, $H^1(\sigma) = {\text{\rm Char}}(G)$. 
\endproclaim
\noindent
{\it Proof}. Theorem 4.1 shows that there 
exist $u\in \Cal H$
such that $w_g=\lambda_g u\sigma_g(u^*)$. On  
the other hand, 
by Remark 2.1, $\lambda_g$ is a character of $G$.

Moreover, if two characters $\lambda_g, \lambda'_g$ are equivalent 
then there exists a unitary element $u \in A$ such that 
$\lambda_g 1= \lambda'_g u\sigma_g(u^*), \forall g\in G$. 
Thus, $\sigma_g(u) \in \Bbb Cu, \forall g\in G$. 
But since $\sigma$ is weakly mixing, the only 
finite dimensional $\sigma$-invariant 
subspace of $A$ is $\Bbb C1$, implying that $u \in \Bbb C1$ 
and $\lambda_g = \lambda'_g$. 
\hfill $\square$

\proclaim{4.3. Corollary} 
The first cohomology group $H^1(\sigma)$ 
of a Bernoulli shift action $\sigma$ of 
$SL(n,\Bbb Z), n\ge 3$ is trivial. More generally, if 
$\Gamma$ is any abelian group, $G=SL(n, \Bbb Z) \times \Gamma$  
and $\sigma$ is a $G$-Bernoulli shift, then $H^1(\sigma) = \hat{\Gamma}$. 
\endproclaim
\noindent
{\it Proof.} Indeed for $n\ge 3$, $SL(n,\Bbb Z)$ 
has the property $T$ of Kazhdan ([K]), 
and by the Nielsen Magnum theorem (see for instance [St]), 
for $n\ge 3$ 
the commutator subgroup of 
$G=SL(n,\Bbb Z)\times \Gamma$ 
is equal to  
$SL(n,\Bbb Z)$. 
Thus the group of characters of $G$ is equal to $\hat{\Gamma}$.
\hfill $\square$

\head  References\endhead

\item{[BMe]} B. Bekka, M. Meyer: 
"Ergodic Theory and Topological dynamics of group actions on
Homogeneous Spaces", London Math Soc Lect. Notes {\bf 269},  
Cambridge University Press, 2000. 

\item{[CW]} A. Connes, B. Weiss: 
{\it Property $(\text{\rm T})$ and
asymptotically 
invariant sequences}, Israel. J. Math. {\bf 37} (1980), 209-210.

\item{[dHV]} P. de la Harpe, A. Valette: ``La propri\'et\'e T 
de Kazhdan pour les 
groupes localement compacts'', Ast\'erisque {\bf 175}, Soc. Math. de France (1989).

\item{[Jo]} P. Jolissaint: {\it On the relative property T}, 
preprint 2001.

\item{[K]} D. Kazhdan: {\it Connection of the dual space of a group 
with the structure of its closed subgroups}, Funct. Anal. and its Appl. 
{\bf1} (1967), 63-65.

\item{[M]} G. Margulis: {\it Finitely-additive invariant measures 
on Euclidian spaces}, Ergodic. Th. and Dynam. Sys. {\bf 2} (1982), 
383-396. 

\item{[MvN]} F. Murray, J. von Neumann:
{\it Rings of operators IV}, Ann. Math. {\bf 44}
(1943), 716-808.

\item{[Po]} S. Popa: {\it Some rigidity 
results for non-commutative Bernoulli shifts}, MSRI preprint, 2001-005.

\item{[S1]} K. Schmidt: {\it Asymptotically invariant 
sequences and an action of $SL(2, \Bbb Z)$ on the 
$2$-sphere}, Israel. J. Math. {\bf 37} (1980), 193-208.

\item{[S2]} K. Schmidt: {\it Amenabilty, Kazhdan's property T, 
strong ergodicity and invariant means 
for ergodic group-actions}, Ergod. Th. \& Dynam. Sys. {\bf 1} 
(1981), 223-236.

\item{[St]} R. Steinberg: {\it Some consequences of 
elementary relations of $SL(n)$},
 Contemporary Math.,{\bf 45} (1985), 335-350.

\enddocument